\documentclass{amsart}

\usepackage{url,amssymb}
\newtheorem{thm}{Theorem}[section]
\newtheorem{prop}[thm]{Proposition}
\newtheorem{lem}[thm]{Lemma}
\theoremstyle{definition}
\newtheorem{defn}[thm]{Definition}
\newtheorem*{fact}{Fact}

\theoremstyle{remark}
\newtheorem{rem}[thm]{Remark}
\numberwithin{equation}{section}

\begin{document}
\title[Finitely presentable Kazhdan groups]
{Finitely presentable, non-Hopfian groups with Kazhdan's
Property~(T) and infinite outer automorphism group}

\author[Cornulier]{Yves de Cornulier}
\date{Oct 28, 2005; Theorem 3.1 corrected, May 12, 2010}
\address{ \'Ecole Polytechnique F\'ed\'erale de Lausanne (EPFL)\\
Institut de G\'eom\'etrie, Alg\`ebre et Topologie (IGAT)\\
CH-1015 Lausanne, Switzerland}

\email{decornul@clipper.ens.fr}

%20G25 Linear algebraic groups over local fields and their integers
%20F28 Automorphism groups of groups
%22E50 Representations of Lie and linear algebraic groups over local fields
%17B56 Cohomology of Lie (super)algebras

\tolerance=5000
% because of some overflows, it is need in order to make use of \linebreak[1]
% (otherwise the dangerous \linebreak is needed) - \allowbreak makes
% the previous line not aligned on the right so I do not use it.

%\commby{Dan M. Barbasch}

\begin{abstract} We give simple examples of Kazhdan groups with
infinite outer automorphism groups. This answers a question of
Paulin, independently answered by Ollivier and Wise by completely
different methods. As arithmetic lattices in (non-semisimple) Lie
groups, our examples are in addition finitely presented.

We also use results of Abels about compact presentability of
$p$-adic groups to exhibit a finitely presented non-Hopfian
Kazhdan group. This answers a question of Ollivier and Wise.
\end{abstract}

\maketitle

%%%%%%%%%%%%%%%%%%%%%%%%%%%%%%%%%%%%%%%%%%%%%%%%%%%%%%%%%%%%%%%%%%%%%%%%%%%%%%%%%%%%%%%%%%%%%%%%%%%%%%
%%%%%%%%%%%%%%%%%%%%%%%%%%%%%%%%%%%%%%%%%%%%%%%%%%%%%%%%%%%%%%%%%%%%%%%%%%%%%%%%%%%%%%%%%%%%%%%%%DOC%%
%%%%%%%%%%%%%%%%%%%%%%%%%%%%%%%%%%%%%%%%%%%%%%%%%%%%%%%%%%%%%%%%%%%%%%%%%%%%%%%%%%%%%%%%%%%%%%%%%%V%%%

\section{Introduction}

Recall that a locally compact group is said to have Property~(T)
if every weakly continuous unitary representation with almost
invariant vectors\footnote{A representation $\pi:G\to
\mathcal{U}(\mathcal{H})$ almost has invariant vectors if for
every $\varepsilon>0$ and every finite subset $F\subseteq G$,
there exists a unit vector $\xi\in\mathcal{H}$ such that
$\|\pi(g)\xi-\xi\|<\varepsilon$ for every $g\in F$.} has nonzero
invariant vectors.

It was asked by Paulin in \cite[p.134]{HV} (1989) whether there
exists a group with Kazhdan's Property~(T) and with infinite outer
automorphism group. This question remained unanswered until 2004;
in particular, it is Question 18 in \cite{W}.

This question was motivated by the two following special cases.
The first is the case of lattices in {\it semisimple} groups over
local fields, which have long been considered as prototypical
examples of groups with Property~(T). If $\Gamma$ is such a
lattice, Mostow's rigidity Theorem and the fact that semisimple
groups have finite outer automorphism group imply that
$\textnormal{Out}(\Gamma)$ is finite. Secondly, a new source of
groups with Property~(T) appeared when Zuk \cite{Zuk} proved that
certain models of random groups have Property~(T). But they are
also hyperbolic, and Paulin proved \cite{Pau} that a hyperbolic
group with Property~(T) has a finite outer automorphism group.

However, it turns out that various arithmetic lattices in
appropriate {\em non-semisimple} groups provide examples. For
instance, consider the additive group
$\textnormal{Mat}_{m\times n}(\mathbf{Z})$ of $m\times n$ matrices over
$\mathbf{Z}$, endowed with the action of
$\textnormal{GL}_n(\mathbf{Z})$ by left multiplication.

\begin{prop}
For every $n\ge 3$, $m\ge 1$,
$\textnormal{SL}_n(\mathbf{Z})\ltimes
\textnormal{Mat}_{m\times n}(\mathbf{Z})$ is a finitely presented linear
group, has Property~(T), is non-coHopfian\footnote{A group is
coHopfian (resp.\ Hopfian) if it is isomorphic to no proper
subgroup (resp.\ quotient) of itself.}, and its outer automorphism
group contains a copy of $\textnormal{PGL}_m(\mathbf{Z})$, hence
is infinite if $m\ge 2$.\label{p:Out_infini}
\end{prop}

We later learned that Ollivier and Wise \cite{OW} had
independently found examples of a very different nature. They
embed any countable group $G$ in $\text{Out}(\Gamma)$, where
$\Gamma$ has Property~(T), is a subgroup of a torsion-free
hyperbolic group, satisfying a certain ``graphical" small
cancelation condition (see also \cite{BS}). In contrast to our
examples, theirs are not, a priori, finitely presented; on the
other hand, our examples are certainly not subgroups of hyperbolic
groups since they all contain a copy of~$\mathbf{Z}^2$.

They also construct in \cite{OW} a non-coHopfian group with
Property~(T) that embeds in a hyperbolic group. Proposition
\ref{p:Out_infini} actually answers two questions in their paper:
namely, whether there exists a finitely presented group with
Property~(T) and without the coHopfian Property (resp.\ with
infinite outer automorphism group).

\begin{rem}Another example of a non-coHopfian group with Property~(T) is
\linebreak[1]$\textnormal{PGL}_n(\mathbf{F}_p[X])$ when $n\ge 3$.
This group is finitely presentable if $n\ge 4$ \cite{RS} (but not
for $n=3$ \cite{Behr}). In contrast with the previous examples,
the Frobenius morphism $\textnormal{Fr}$ induces an isomorphism
onto a subgroup of {\em infinite} index, and the intersection
$\bigcap_{k\ge 0}\textnormal{Im}(\textnormal{Fr}^k)$ is reduced
to~$\{1\}$.
\end{rem}

Ollivier and Wise also constructed in \cite{OW} the first examples
of non-Hopfian groups with Property~(T). They asked whether a
finitely presented example exists. Although linear finitely
generated groups are residually finite, hence Hopfian, we use them
to positively answer their question.

\begin{thm}
There exists a $S$-arithmetic lattice $\Gamma$, and a central
subgroup $Z\subset \Gamma$, such that $\Gamma$ and $\Gamma/Z$ are
finitely presented, have Property~(T), and $\Gamma/Z$ is
non-Hopfian.\label{t:nonhopf}
\end{thm}

The group $\Gamma$ has a simple description as a matrix group from
which Property~(T) and the non-Hopfian property for $\Gamma/Z$ are
easily checked (Proposition \ref{p:hopfT}). Section \ref{s:3} is
devoted to prove finite presentability of $\Gamma$. We use here a
general criterion for finite presentability of $S$-arithmetic
groups, due to Abels \cite{Abels}. It involves the computation of
the first and second cohomology group of a suitable Lie algebra.

\section{Proofs of all results except finite presentability of~$\Gamma$}

We need some facts about Property~(T).

\begin{lem}[see {\cite[Chap. 3, Th\'eor\`eme 4]{HV}}]
Let $G$ be a locally compact group, and $\Gamma$ a lattice in $G$.
Then $G$ has Property~(T) if and only if $\Gamma$ has
Property~(T).\qed\label{inherited_lattices}
\end{lem}

The next lemma is an immediate consequence of the classification
of semisimple algebraic groups over local fields with Property (T)
(see \cite[Chap.~III, Theorem~5.6]{Margulis}) and S.~P.~Wang's
results on the non-semisimple case \cite[Theorem 2.10]{Wang}.

\begin{lem}
Let $\mathbf{K}$ be a local field, $G$ a connected linear
algebraic group defined over $\mathbf{K}$. Suppose that $G$ is
perfect, and, for every simple quotient $S$ of $G$, either $S$ has
$\mathbf{K}$-rank $\ge 2$, or $\mathbf{K}=\mathbf{R}$ and $S$ is
isogeneous to either $\textnormal{Sp}(n,1)$ ($n\ge 2$) or
$\textnormal{F}_{4(-20)}$. If $\textnormal{char}(\mathbf{K})>0$,
suppose in addition that $G$ has a Levi decomposition defined over
$\mathbf{K}$. Then $G(\mathbf{K})$ has
Property~(T).\qed\label{wang_perfect}
\end{lem}

\begin{proof}[Proof of Proposition \ref{p:Out_infini}] The group
$\textnormal{SL}_n(\mathbf{Z})\ltimes
\textnormal{Mat}_{m\times n}(\mathbf{Z})$ is linear in dimension $n+m$.
As a semidirect product of two finitely presented groups, it is
finitely presented. For every $k\ge 2$, it is isomorphic to its
proper subgroup $\textnormal{SL}_n(\mathbf{Z})\ltimes
k\textnormal{Mat}_{m\times n}(\mathbf{Z})$ of finite index $k^{mn}$.

The group $\textnormal{GL}_m(\mathbf{Z})$ acts on
$\textnormal{Mat}_{m\times n}(\mathbf{Z})$ by right multiplication. Since
this action commutes with the left multiplication of
$\textnormal{SL}_n(\mathbf{Z})$, $\textnormal{GL}_m(\mathbf{Z})$
acts on the semidirect product
$\textnormal{SL}_n(\mathbf{Z})\ltimes
\textnormal{Mat}_{m\times n}(\mathbf{Z})$ by automorphisms, and, by an
immediate verification, this gives an embedding of
$\textnormal{GL}_m(\mathbf{Z})$ if $n$ is odd or
$\textnormal{PGL}_m(\mathbf{Z})$ if $n$ is even into
$\textnormal{Out}(\textnormal{SL}_n(\mathbf{Z})\ltimes
\textnormal{Mat}_{m\times n}(\mathbf{Z}))$ (it can be shown that this is
an isomorphism if $n$ is odd; if $n$ is even, the image has index
two). In particular, if $m\ge 2$, then
$\textnormal{SL}_n(\mathbf{Z})\ltimes
\textnormal{Mat}_{m\times n}(\mathbf{Z})$ has infinite outer automorphism
group.

On the other hand, in view of Lemma \ref{inherited_lattices}, it
has Property~(T) (actually for all $m\ge 0$): indeed,
$\textnormal{SL}_n(\mathbf{Z})\ltimes
\textnormal{Mat}_{m\times n}(\mathbf{Z})$ is a lattice in
$\textnormal{SL}_n(\mathbf{R})\ltimes
\textnormal{Mat}_{m\times n}(\mathbf{R})$, which has Property~(T) by
Lemma \ref{wang_perfect} as $n\ge 3$.\end{proof}

We now turn to the proof of Theorem \ref{t:nonhopf}. The following
lemma is immediate, and already used
in \cite[Th.~4(iii)]{Hall} and~\cite{AbelsPapier}.%REF HALL?

\begin{lem}
Let $\Gamma$ be a group, $Z$ a central subgroup. Let $\alpha$ be
an automorphism of $\Gamma$ such that $\alpha(Z)$ is a proper
subgroup of $Z$. Then $\alpha$ induces a surjective, non-injective
endomorphism of $\Gamma/Z$, whose kernel is
$\alpha^{-1}(Z)/Z$.\qed\label{l:nonhopf}
\end{lem}

\begin{defn}
Fix $n_1,n_2,n_3,n_4\in\mathbf{N}-\{0\}$ with $n_2,n_3\ge 3$. We
set $\Gamma=G(\mathbf{Z}[1/p])$, where $p$ is any prime, and $G$
is algebraic the group defined as matrices by blocks of size
$n_1,n_2,n_3,n_4$:

$$\begin{pmatrix}
  I_{n_1} & (*)_{12} & (*)_{13} & (*)_{14} \\
  0 & (**)_{22} & (*)_{23} & (*)_{24} \\
  0 & 0 & (**)_{33} & (*)_{34} \\
  0 & 0 & 0 & I_{n_4} \\
\end{pmatrix},$$
where $(*)$ denote any matrices and $(**)_{ii}$ denote matrices in
$\textnormal{SL}_{n_i}$, $i=2,3$.

The centre of $G$ consists of matrices of the form
$\begin{pmatrix}
  I_{n_1} & 0 & 0 & (*)_{14} \\
  0 & I_{n_2} & 0 & 0 \\
  0 & 0 & I_{n_3} & 0 \\
  0 & 0 & 0 & I_{n_4} \\
\end{pmatrix}$. Define $Z$
as the centre of $G(\mathbf{Z})$.\label{d:Gamma_nonhopfien}
\end{defn}

\begin{rem}This group is related to an example of Abels: in
\cite{AbelsPapier} he considers the same group, but with blocks
$1\times 1$, and $\textnormal{GL}_1$ instead of
$\textnormal{SL}_1$ in the diagonal. Taking the points over
$\mathbf{Z}[1/p]$, and taking the quotient by a cyclic subgroup if
the centre, this provided the first example of a finitely
presentable non-Hopfian solvable group.\end{rem}

\begin{rem}
If we do not care about finite presentability, we can take $n_3=0$ (i.e. 3 blocks suffice), as in P.~Hall's original solvable example \cite[Th.~4(iii)]{Hall}.
\end{rem}

We begin by easy observations. Identify $\textnormal{GL}_{n_1}$ to
the upper left diagonal block. It acts by \textit{conjugation} on
$G$ as follows:
$$\begin{pmatrix}
  u & 0 & 0 & 0 \\
  0 & I & 0 & 0 \\
  0 & 0 & I & 0 \\
  0 & 0 & 0 & I \\
\end{pmatrix}\cdot\begin{pmatrix}
  I & A_{12} & A_{13} & A_{14} \\
  0 & B_{2} & A_{23} & A_{24} \\
  0 & 0 & B_3 & A_{34} \\
  0 & 0 & 0 & I \\
\end{pmatrix}=
\begin{pmatrix}
  I & uA_{12} & uA_{13} & uA_{14} \\
  0 & B_{2} & A_{23} & A_{24} \\
  0 & 0 & B_3 & A_{34} \\
  0 & 0 & 0 & I \\
\end{pmatrix}.$$

This gives an action of $\textnormal{GL}_{n_1}$ on $G$, and also
on its centre, and this latter action is faithful. In particular,
for every commutative ring $R$, $\textnormal{GL}_{n_1}(R)$ embeds
in $\textnormal{Out}(G(R))$.

From now on, we suppose that $R=\mathbf{Z}[1/p]$, and
$u=pI_{n_1}$. The automorphism of $\Gamma=G(\mathbf{Z}[1/p])$
induced by $u$ maps $Z$ to its proper subgroup $Z^p$. In view of
Lemma \ref{l:nonhopf}, this implies that $\Gamma/Z$ is
non-Hopfian.

\begin{prop}
The groups $\Gamma$ and $\Gamma/Z$ are finitely generated, have
Property~(T), and $\Gamma/Z$ is non-Hopfian.\label{p:hopfT}
\end{prop}
\begin{proof} We have just proved that $\Gamma/Z$ is non-Hopfian. By the
Borel-Harish-Chandra Theorem \cite{BHC}, $\Gamma$ is a lattice in
$G(\mathbf{R})\times G(\mathbf{Q}_p)$. This group has Property~(T)
as a consequence of Lemma \ref{wang_perfect}. By Lemma
\ref{inherited_lattices}, $\Gamma$ also has Property~(T). Finite
generation is a consequence of Property~(T) \cite[Lemme 10]{HV}.
Since Property~(T) is (trivially) inherited by quotients,
$\Gamma/Z$ also has Property~(T).\end{proof}

\begin{rem}
This group has a surjective endomorphism with nontrivial finite
kernel. We have no analogous example with infinite kernel. Such
examples might be constructed if we could prove that some groups
over rings of dimension $\ge 2$ such as
$\textnormal{SL}_n(\mathbf{Z}[X])$ or
$\textnormal{SL}_n(\mathbf{F}_p[X,Y])$ have Property~(T), but this
is an open problem \cite{Sha}. The non-Hopfian Kazhdan group of
Ollivier and Wise \cite{OW} is torsion-free, so the kernel is
infinite in their case.
\end{rem}

\begin{rem}
It is easy to check that
$\textnormal{GL}_{n_1}(\mathbf{Z})\times\textnormal{GL}_{n_4}(\mathbf{Z})$
embeds in $\textnormal{Out}(\Gamma)$ and
$\textnormal{Out}(\Gamma/Z)$. In particular, if $\max(n_1,n_4)\ge
2$, then these outer automorphism groups are infinite.
\end{rem}

We finish this section by observing that $Z$ is a finitely
generated subgroup of the centre of $\Gamma$, so that finite
presentability of $\Gamma/Z$ immediately follows from that
of~$\Gamma$.

\section{Finite presentability of $\Gamma$}\label{s:3}

We recall that a Hausdorff topological group $H$ is
\textit{compactly presented} if there exists a compact generating
subset $C$ of $H$ such that the abstract group $H$ is the quotient
of the group freely generated by $C$ by relations of bounded
length. See Abels \cite[\S 1.1]{Abels} for more about compact
presentability.

Kneser \cite{Kneser} has proved that for every linear algebraic
$\mathbf{Q}_p$-group, the $S$-arithmetic lattice
$G(\mathbf{Z}[1/p])$ is finitely presented if and only if
$G(\mathbf{Q}_p)$ is compactly presented. A characterization of
the linear algebraic $\mathbf{Q}_p$-groups $G$ such that
$G(\mathbf{Q}_p)$ is compactly presented is given in~\cite{Abels}.
This criterion requires the study of a solvable cocompact subgroup
of $G(\mathbf{Q}_p)$, which seems tedious to carry out in our
specific example.

Let us describe another sufficient criterion for compact
presentability, also given in \cite{Abels}, which is applicable to
our example. Let $U$ be the unipotent radical in $G$, and let $S$
denote a Levi factor defined over $\mathbf{Q}_p$, so that
$G=S\ltimes U$. Let $\mathfrak{u}$ be the Lie algebra of $U$, and
$D$ be a maximal $\mathbf{Q}_p$-split torus in $S$. We recall that
the first homology group of $\mathfrak{u}$ is defined as the
abelianization
$$H_1(\mathfrak{u})=\mathfrak{u}/[\mathfrak{u},\mathfrak{u}],$$
and the second homology group of $\mathfrak{u}$ is defined as
$\textnormal{Ker}(d_2)/\textnormal{Im}(d_3)$, where the maps
$$\mathfrak{u}\wedge\mathfrak{u}\wedge\mathfrak{u}\stackrel{d_3}
\to\mathfrak{u}\wedge\mathfrak{u}\stackrel{d_2}\to\mathfrak{u}$$
are defined by:

$$d_2(x_1\wedge x_2)=-[x_1,x_2]\;\;\text{and}\;\; d_3(x_1\wedge x_2\wedge
x_3)=x_3 \wedge [x_1,x_2]+x_2\wedge [x_3,x_1]+x_1\wedge
[x_2,x_3].$$

We can now state the result by Abels that we use (see
\cite[Theorem 6.4.3 and Remark 6.4.5]{Abels}).

\begin{thm}[Abels]\label{abth}\footnote{In the published version, the corresponding theorem is a misquotation of Abels' theorem. This version is corrected. The numbering of the lemmas is not affected. The reference to the erratum is published as Proc. Amer. Math. Soc. 139 (2011), 383--384. It essentially includes the correct statement of Theorem \ref{abth} (as given here) and the subsequent lines ``on the proof". The verifications in the published version were enough to be applicable to the correct version of Theorem \ref{abth}, so no change in the computations has been done.}

Let $G$ be a connected
linear algebraic group over~$\mathbf{Q}_p$. Suppose that $G$ is unipotent-by-semisimple, i.e.\ $G=US$, where $U$ is the unipotent radical and $S$ is a semisimple Levi factor. We furthermore assume that $S$ is split semisimple without any simple factors of rank one.
Then $G(\mathbf{Q}_p)$ is compactly presented if and only if the following two conditions are satisfied
\begin{itemize}

\item[(i)] $H_1(\mathfrak{u})^S=\{0\}$;

\item[(ii)] $H_2(\mathfrak{u})^S=\{0\}$.
\end{itemize}
\label{t:Abels}\end{thm}
\begin{proof}[On the proof]
This relies on \cite[Theorem 6.4.3 and Remark 6.4.5]{Abels}. A few comments are necessary: 
\begin{itemize}
\item Condition (1a) of \cite[Remark 6.4.5]{Abels} involves the orthogonal $\Phi^\bot$ of the subspace generated by roots. It states that $\Phi^\bot$ does not contain dominant weights $\omega_1,\omega_2$ of the $S$-module $H_1(\mathfrak{u})$ with $0\in [\omega_1,\omega_2]$. Since $S$ is semisimple, $\Phi^\bot=\{0\}$ and (1a) just means that 0 is not a dominant weight, which is exactly Condition (i) above.
(Note that (i) is actually a necessary and sufficient condition for $G(\mathbf{Q}_p)$ to be compactly {\em generated}, see \cite[Theorem 6.4.4]{Abels}.)
\item As noticed in \cite[p.~132]{Abels}, Condition~(1b) of \cite[Remark 6.4.5]{Abels} is superfluous when $S$ has no factors of rank $\le 1$.
\item (ii) is a restatement of Condition~2 of \cite[Theorem 6.4.3]{Abels}.\qedhere
\end{itemize}
\end{proof}

%\begin{thm} Let $G$ be a connected
%linear algebraic group over~$\mathbf{Q}_p$. Suppose that the
%following assumptions are fulfilled:

%\begin{itemize}

%\item[(i)] $G$ is $\mathbf{Q}_p$-split.

%\item[(ii)] $G$ has no simple quotient of $\mathbf{Q}_p$-rank one.

%\item[(iii)] 0 does not lie on the segment joining two dominant
%weights for the adjoint representation of $S$ on
%$H_1(\mathfrak{u})$.

%\item[(iv)] 0 is not a dominant weight for an irreducible
%subrepresentation of the adjoint representation of $S$ on
%$H_2(\mathfrak{u})$.\end{itemize}

%Then $G(\mathbf{Q}_p)$ is compactly presented.\qed\label{t:Abels}
%\end{thm}

We now return to our particular example $G$ from Definition \ref{d:Gamma_nonhopfien}; it is unipotent-by-semisimple and the semisimple Levi factor $S=\textnormal{SL}_{n_2}\times\textnormal{SL}_{n_3}$ is split with no factor of rank one, so it fulfills the assumptions of Theorem \ref{abth}. So its compact presentability is equivalent to Conditions (i) and (ii) of this theorem. Keep the previous notation $S$, $D$, $U$,
$\mathfrak{u}$, so that $S$ (resp.\ $D$) denoting in our case the
diagonal by blocks (resp.\ diagonal) matrices in $G$, and $U$
denotes the matrices in $G$ all of whose diagonal blocks are the
identity. The set of indices of the matrix is partitioned as
$I=I_1\sqcup I_2\sqcup I_3\sqcup I_4$, with $|I_j|=n_j$ as in
Definition \ref{d:Gamma_nonhopfien}. It follows that, for every
field $K$,

$$\mathfrak{u}(K)=\left\{T\in \text{End}(K^I),\;\forall j,\;T(K^{I_j})\subset
\bigoplus_{i<j}K^{I_i}\right\}.$$

Throughout, we use the following notation: a letter such as $i_k$
(or $j_k$, etc.) implicitly means $i_k\in I_k$. Define, in an
obvious way, subgroups $U_{ij}$, $i<j$, of $U$, and their Lie
algebras~$\mathfrak{u}_{ij}$.

%%%%%%%%%%%%%

We begin by checking Condition (i) of Theorem \ref{t:Abels}. This follows from the following lemma.

\begin{lem}\footnote{Compared to the published version, only the second sentence in the statement of Lemma \ref{l:H1_criterion} has been added, and the proof has not been modified, although the statement is stronger than what is actually needed for the corrected version given here of Theorem \ref{t:Abels}.}
For any two weights of the action of $D$ on $H_1(\mathfrak{u})$, 0
is not on the segment joining them. In particular, 0 is not a weight of the action of $D$ on $H_1(\mathfrak{u})$.
\label{l:H1_criterion}
\end{lem}
\begin{proof} Recall that $H_1(\mathfrak{u})=\mathfrak{u}/[\mathfrak{u},\mathfrak{u}]$. So it
suffices to look at the action on the supplement $D$-subspace
$\mathfrak{u}_{12}\oplus\mathfrak{u}_{23}\oplus\mathfrak{u}_{34}$
of $[\mathfrak{u},\mathfrak{u}]$. Identifying $S$ with
$\textnormal{SL}_{n_2}\times\textnormal{SL}_{n_3}$, we denote
$(A,B)$ an element of $D\subset S$. We also denote by $e_{pq}$ the
matrix whose coefficient $(p,q)$ equals one and all others are
zero.

$$(A,B)\cdot e_{i_1j_2}=a_{j_2}^{-1}e_{i_1j_2},\quad(A,B)\cdot
e_{j_2k_3}=a_{j_2}b_{k_3}^{-1}e_{j_2k_3},\quad(A,B)\cdot
e_{k_3\ell_4}=b_{k_3}e_{k_3\ell_4}.$$

Since $S=\textnormal{SL}_{n_2}\times \textnormal{SL}_{n_3}$, the
weights for the adjoint action on
$\mathfrak{u}_{12}\oplus\mathfrak{u}_{23}\oplus\mathfrak{u}_{34}$
live in $M/P$, where $M$ is the free $\mathbf{Z}$-module of rank
$n_2+n_3$ with basis $(u_1,\dots,u_{n_2},v_1,\dots,v_{n_3})$, and
$P$ is the plane generated by $\sum_{j_2} u_{j_2}$ and
$\sum_{k_3}v_{k_3}$. Thus, the weights are (modulo $P$)
$-u_{j_2}$, $u_{j_2}-v_{k_3}$, $v_{k_3}$ ($1\le j_2\le n_2$, $1\le
k_3\le n_3$).

Using that $n_2,n_3\ge 3$, it is clear that no nontrivial positive
combination of two weights (viewed as elements of
$\mathbf{Z}^{n_2+n_3}$) lies in $P$.\end{proof}

We must now check Condition (ii) of Theorem \ref{t:Abels}, and
therefore compute $H_2(\mathfrak{u})$ as a $D$-module.

\begin{lem}
$\textnormal{Ker}(d_2)\subset\mathfrak{u}\wedge\mathfrak{u}$ is linearly spanned by

\begin{itemize}

\item[(1)] $\mathfrak{u}_{12}\wedge \mathfrak{u}_{12}$,
$\mathfrak{u}_{23}\wedge \mathfrak{u}_{23}$,
$\mathfrak{u}_{34}\wedge \mathfrak{u}_{34}$,
$\mathfrak{u}_{13}\wedge\mathfrak{u}_{23}$,
$\mathfrak{u}_{23}\wedge\mathfrak{u}_{24}$,
$\mathfrak{u}_{12}\wedge\mathfrak{u}_{13}$,
$\mathfrak{u}_{24}\wedge\mathfrak{u}_{34}$,
$\mathfrak{u}_{12}\wedge\mathfrak{u}_{34}$.

\item[(2)] $\mathfrak{u}_{14}\wedge\mathfrak{u}$,
$\mathfrak{u}_{13}\wedge \mathfrak{u}_{13}$,
$\mathfrak{u}_{24}\wedge \mathfrak{u}_{24}$,
$\mathfrak{u}_{13}\wedge \mathfrak{u}_{24}$.

\item[(3)] $e_{i_1j_2}\wedge e_{k_2\ell_3}$ ($j_2\neq k_2$),
$e_{i_2j_3}\wedge e_{k_3\ell_4}$ ($j_3\neq \ell_3$).

\item[(4)] $e_{i_1j_2}\wedge e_{k_2\ell_4}$ ($j_2\neq k_2$),
$e_{i_1j_3}\wedge e_{k_3\ell_4}$ ($j_3\neq k_3$).

\item[(5)] Elements of the form
$\sum_{j_2}\alpha_{j_2}(e_{i_1j_2}\wedge e_{j_2k_3})$ if
$\sum_{j_2}\alpha_{j_2}=0$, and
$\sum_{j_3}\alpha_{j_3}(e_{i_2j_3}\wedge e_{j_3k_4})$ if
\linebreak[1]$\sum_{j_3}\alpha_{j_3}=\nolinebreak 0$.

\item[(6)] Elements of the form
$\sum_{j_2}\alpha_{j_2}(e_{i_1j_2}\wedge
e_{j_2k_4})+\sum_{j_3}\beta_{j_3}(e_{i_1j_3}\wedge e_{j_3k_4})$ if
$\sum_{j_2}\alpha_{j_2}+\sum_{j_3}\beta_{j_3}=0$.
\end{itemize}
\label{kerd2}
\end{lem}

\begin{proof} First observe that $\textnormal{Ker}(d_2)$ contains
$\mathfrak{u}_{ij}\wedge\mathfrak{u}_{kl}$ when
$[\mathfrak{u}_{ij},\mathfrak{u}_{kl}]=0$. This corresponds to (1)
and (2). The remaining cases are
$\mathfrak{u}_{12}\wedge\mathfrak{u}_{23}$,
$\mathfrak{u}_{23}\wedge\mathfrak{u}_{34}$,
$\mathfrak{u}_{12}\wedge\mathfrak{u}_{24}$,
$\mathfrak{u}_{13}\wedge\mathfrak{u}_{34}$.

On the one hand, $\textnormal{Ker}(d_2)$ also contains
$e_{i_1j_2}\wedge e_{k_2\ell_3}$ if $j_2\neq k_2$, etc.; this
corresponds to elements in (3), (4). On the other hand,
$d_2(e_{i_1j_2}\wedge e_{j_2k_3})=-e_{i_1k_3}$,
$d_2(e_{i_2j_3}\wedge e_{j_3k_4})=-e_{i_2k_4}$,
$d_2(e_{i_1j_2}\wedge e_{j_2k_4})=-e_{i_1k_4}$,
$d_2(e_{i_1j_3}\wedge e_{j_3k_4})=-e_{i_1k_4}$. The lemma
follows.\end{proof}

\begin{defn}
Denote by $\mathfrak{b}$ (resp.\ $\mathfrak{h}$) the subspace
spanned by elements in (2), (4), and (6) (resp.\ in (1), (3), and
(5)) of Lemma \ref{kerd2}.
\end{defn}

\begin{prop}
$\textnormal{Im}(d_3)=\mathfrak{b}$, and
$\textnormal{Ker}(d_2)=\mathfrak{b}\oplus\mathfrak{h}$ as
$D$-module. In particular, $H_2(\mathfrak{u})$ is isomorphic to
$\mathfrak{h}$ as a $D$-module.
\end{prop}
\begin{proof} We first prove, in a series of facts, that
$\textnormal{Im}(d_3)\supset\mathfrak{b}$.

\begin{fact}$\mathfrak{u}_{14}\wedge\mathfrak{u}$ is contained in
$\textnormal{Im}(d_3)$.\end{fact}

\begin{proof} If $z\in \mathfrak{u}_{14}$, then $d_3(x\wedge y\wedge z)=z\wedge
[x,y]$. This already shows that $\mathfrak{u}_{14}\wedge
(\mathfrak{u}_{13}\oplus \mathfrak{u}_{24} \oplus
\mathfrak{u}_{14})$ is contained in $\textnormal{Im}(d_3)$, since
$[\mathfrak{u},\mathfrak{u}]=\mathfrak{u}_{13}\oplus\mathfrak{u}_{24}
\oplus \mathfrak{u}_{14}$.

Now, if $(x,y,z)\in
\mathfrak{u}_{24}\times\mathfrak{u}_{12}\times\mathfrak{u}_{34}$,
then $d_3(x\wedge y\wedge z)=z\wedge [x,y]$. Since
$[\mathfrak{u}_{24},\mathfrak{u}_{12}]=\mathfrak{u}_{14}$, this
implies that $\mathfrak{u}_{14}\wedge\mathfrak{u}_{34}\subset
\textnormal{Im}(d_3)$. Similarly,
$\mathfrak{u}_{14}\wedge\mathfrak{u}_{12}\subset
\textnormal{Im}(d_3)$.

Finally we must prove that
$\mathfrak{u}_{14}\wedge\mathfrak{u}_{23}\subset
\textnormal{Im}(d_3)$. This follows from the formula
$e_{i_1j_4}\wedge e_{k_2\ell_3}=d_3(e_{i_1m_2}\wedge
e_{k_2\ell_3}\wedge e_{m_2j_4})$, where $m_2\neq k_2$ (so that we
use that $|I_2|\ge 2$).\end{proof}

\begin{fact}
$\mathfrak{u}_{13}\wedge\mathfrak{u}_{13}$ and, similarly,
$\mathfrak{u}_{24}\wedge\mathfrak{u}_{24}$, are contained in
$\textnormal{Im}(d_3)$.
\end{fact}
\begin{proof} If $(x,y,z)\in
\mathfrak{u}_{12}\times\mathfrak{u}_{23}\times\mathfrak{u}_{13}$,
then $d_3(x\wedge y\wedge z)=z\wedge [x,y]$. Since
$[\mathfrak{u}_{12},\mathfrak{u}_{23}]=\mathfrak{u}_{13}$, this
implies that $\mathfrak{u}_{13}\wedge\mathfrak{u}_{13}\subset
\textnormal{Im}(d_3)$.\end{proof}

\begin{fact}
$\mathfrak{u}_{13}\wedge\mathfrak{u}_{24}$ is contained in
$\textnormal{Im}(d_3)$.
\end{fact}
\begin{proof} $d_3(e_{i_1k_2}\wedge e_{k_2\ell_3}\wedge
e_{k_2j_4})=e_{k_2j_4}\wedge e_{i_1\ell_3}+e_{i_1j_4}\wedge
e_{k_2\ell_3}$. Since we already know that $e_{i_1j_4}\wedge
e_{k_2\ell_3}\in\textnormal{Im}(d_3)$, this implies
$e_{k_2j_4}\wedge
e_{i_1\ell_3}\in\textnormal{Im}(d_3)$.\end{proof}

\begin{fact}
The elements in (4) are in $\textnormal{Im}(d_3)$.
\end{fact}
\begin{proof} $d_3(e_{i_1j_2}\wedge e_{j_2k_3}\wedge e_{\ell_3m_4})=
-e_{i_1k_3}\wedge e_{\ell_3m_4}$ if $k_3\neq\ell_3$. The other
case is similar.\end{proof}

\begin{fact}
The elements in (6) are in $\textnormal{Im}(d_3)$.
\end{fact}
\begin{proof} $d_3(e_{i_1j_2}\wedge e_{j_2k_3}\wedge e_{k_3\ell_4})=
-e_{i_1k_3}\wedge e_{k_3\ell_4}+e_{i_1j_2}\wedge e_{j_2\ell_4}$.
Such elements linearly span all elements as in (6).\end{proof}

Conversely, we must check
$\textnormal{Im}(d_3)\subset\mathfrak{b}$. By straightforward
verifications:

\begin{itemize}

\item
$d_3(\mathfrak{u}_{14}\wedge\mathfrak{u}\wedge\mathfrak{u})
\subset\mathfrak{u}_{14}\wedge\mathfrak{u}$.

\item
$d_3(\mathfrak{u}_{13}\wedge\mathfrak{u}_{23}\wedge\mathfrak{u}_{24})=0$

\item
$d_3(\mathfrak{u}_{12}\wedge\mathfrak{u}_{13}\wedge\mathfrak{u}_{24})$,
$d_3(\mathfrak{u}_{13}\wedge\mathfrak{u}_{24}\wedge\mathfrak{u}_{34})$,
$d_3(\mathfrak{u}_{12}\wedge\mathfrak{u}_{13}\wedge\mathfrak{u}_{34})$,
$d_3(\mathfrak{u}_{12}\wedge\mathfrak{u}_{24}\wedge\mathfrak{u}_{34})$
are all contained in $\mathfrak{u}_{14}\wedge\mathfrak{u}$.

\item
$d_3(\mathfrak{u}_{12}\wedge\mathfrak{u}_{13}\wedge\mathfrak{u}_{23})\subset
\mathfrak{u}_{13}\wedge\mathfrak{u}_{13}$, and similarly
$d_3(\mathfrak{u}_{23}\wedge\mathfrak{u}_{24}\wedge\mathfrak{u}_{34})\subset
\mathfrak{u}_{24}\wedge\mathfrak{u}_{24}$.

\item
$d_3(\mathfrak{u}_{12}\wedge\mathfrak{u}_{23}\wedge\mathfrak{u}_{24})$
and similarly
$d_3(\mathfrak{u}_{13}\wedge\mathfrak{u}_{23}\wedge\mathfrak{u}_{34})$
are contained in
$\mathfrak{u}_{14}\wedge\mathfrak{u}_{23}+\mathfrak{u}_{13}\wedge\mathfrak{u}_{24}$.

\item The only remaining case is that of
$\mathfrak{u}_{12}\wedge\mathfrak{u}_{23}\wedge\mathfrak{u}_{34}$:
$d_3(e_{i_1j_2}\wedge e_{j'_2k_3}\wedge
e_{k'_3\ell_4})=\delta_{k_3k'_3}e_{i_1j_2}\wedge
e_{j'_2\ell_4}-\delta_{j_2j'_2} e_{i_1k_3}\wedge e_{k'_3\ell_4}$,
which lies in (4) or in (6).

\end{itemize}

Finally $\textnormal{Im}(d_3)=\mathfrak{b}$.

\medskip

It follows from Lemma \ref{kerd2} that
$\textnormal{Ker}(d_2)=\mathfrak{h}\oplus\mathfrak{b}$. Since
$\mathfrak{b}=\textnormal{Im}(d_3)$, this is a $D$-submodule. Let
us check that $\mathfrak{h}$ is also a $D$-submodule; the
computation will be used in the sequel.

The action of $S$ on $\mathfrak{u}$ by {\em conjugation} is given
by:

$$\begin{pmatrix}
  1 & 0 & 0 & 0 \\
  0 & A & 0 & 0 \\
  0 & 0 & B & 0 \\
  0 & 0 & 0 & 1 \\
\end{pmatrix}\cdot\begin{pmatrix}
  0 & X_{12} & X_{13} & X_{14} \\
  0 & 0 & X_{23} & X_{24} \\
  0 & 0 & 0 & X_{34} \\
  0 & 0 & 0 & 0 \\
\end{pmatrix}=\begin{pmatrix}
  0 & X_{12}A^{-1} & X_{13}B^{-1} & X_{14} \\
  0 & 0 & AX_{23}B^{-1} & AX_{24} \\
  0 & 0 & 0 & BX_{34} \\
  0 & 0 & 0 & 0 \\
\end{pmatrix}$$

We must look at the action of $D$ on the elements in (1), (3), and
(5). We fix $(A,B)\in D\subset
S\simeq\textnormal{SL}_{n_2}\times\textnormal{SL}_{n_3}$, and we
write $A=\sum_{j_2}a_{j_2}e_{j_2j_2}$ and
$B=\sum_{k_3}b_{k_3}e_{k_3k_3}$.

\begin{itemize}

\item (1):
\begin{equation}(A,B)\cdot e_{i_1j_2}\wedge
e_{k_1\ell_2}=e_{i_1j_2}A^{-1}\wedge
e_{k_1\ell_2}A^{-1}=a_{j_2}^{-1}a_{\ell_2}^{-1}e_{i_1j_2}\wedge
e_{k_1\ell_2}.\label{e_2}\end{equation} The action on other
elements in (1) has a similar form.

\item (3) ($j_2\neq k_2$):
\begin{equation}(A,B)\cdot e_{i_1j_2}\wedge
e_{k_2\ell_3}=e_{i_1j_2}A^{-1}\wedge
Ae_{k_2\ell_4}B^{-1}=a_{j_2}^{-1}a_{k_2}b_{\ell_3}^{-1}e_{i_1j_2}\wedge
e_{k_2\ell_3}.\label{e_3}\end{equation} The action on the other
elements in (3) has a similar form.

\item (5) ($\sum_{j_2}\alpha_{j_2}=0$)

\begin{align}\nonumber (A,B)\cdot\sum_{j_2}\alpha_{j_2}(e_{i_1j_2}\wedge
e_{j_2k_3})\;\;=\;\;\sum_{j_2}\alpha_{j_2}(e_{i_1j_2}A^{-1}\wedge
Ae_{j_2k_3}B^{-1})\\=\sum_{j_2}\alpha_{j_2}a_{j_2}^{-1}(e_{i_1j_2}\wedge
a_{j_2}b_{k_3}^{-1}e_{j_2k_3})\;\;=\;\;b_{k_3}^{-1}\left(\sum_{j_2}\alpha_{j_2}(e_{i_1j_2}\wedge
e_{j_2k_3})\right).\label{e_5}\end{align}
   The other case in (5) has a similar form.\qedhere\end{itemize}\end{proof}

\begin{lem}
0 is not a weight for the action of $D$ on
$H_2(\mathfrak{u})$.\label{l:H2_criterion}
\end{lem}
\begin{proof} As described in the proof of Lemma \ref{l:H1_criterion}, we
think of weights as elements of $M/P$. Hence, we describe weights
as elements of $M=\mathbf{Z}^{n_2+n_3}$ rather than $M/P$, and
must check that no weight lies in $P$.

\begin{itemize}

\item[(1)] In (\ref{e_2}), the weight is $-u_{j_2}-u_{\ell_2}$,
hence does not belong to $P$ since $n_2\ge 3$. The other
verifications are similar.

\item[(3)] In (\ref{e_3}), the weight is
$-u_{j_2}+u_{k_2}-v_{\ell_3}$, hence does not belong to $P$. The
other verification for (3) is similar.

\item[(5)] In (\ref{e_5}), the weight is $-v_{k_3}$, hence does
dot belong to $P$. The other verification is
similar.\qedhere\end{itemize}
\end{proof}

%%%%%%%%%%%%%%

Finally, Lemmas \ref{l:H1_criterion} and \ref{l:H2_criterion}
imply that the conditions (i) and (ii) of Theorem \ref{t:Abels} are satisfied,
so that $\Gamma$ is finitely presented.\qed

%%%%%%%%%%%%%%%%%%%%%%%%%%%%%%%%%%%%%%%%%%%%%%%%%%%%%%%%%%%%%%%%%%%%%%%%%%%%%%%%%%%%%%%%%%%%

\bigskip

\noindent {\bf Acknowledgments.} I thank Herbert Abels, Yann
Ollivier, and Fr\'ed\'eric Paulin for useful discussions, and
Laurent Bartholdi and the referee for valuable comments and
corrections.

\bibliographystyle{amsplain}

%----------------------------------------------------------------------- % End of proc-l.template
%-----------------------------------------------------------------------
\end{document}